\documentclass[11pt,reqno,tbtags]{amsart}
\usepackage{amssymb}
\numberwithin{equation}{section}

\newtheorem{theorem}{Theorem}[section]
\newtheorem{lemma}[theorem]{Lemma}
\newtheorem{proposition}[theorem]{Proposition}

\theoremstyle{definition}

\newtheorem{remark}[theorem]{Remark}

\newtheorem*{ack}{Acknowledgement}

\theoremstyle{remark}

\newenvironment{romenumerate}{\begin{enumerate}
 }{\end{enumerate}}

\newcounter{thmenumerate}
\newenvironment{thmenumerate}
{\setcounter{thmenumerate}{0}%
 \def\item{\par
 \refstepcounter{thmenumerate}\textup{(\roman{thmenumerate})\enspace}}
}
{}

\newcounter{xenumerate}   

\newcommand{\refT}[1]{Theorem~\ref{#1}}

\newcommand{\refL}[1]{Lemma~\ref{#1}}
\newcommand{\refR}[1]{Remark~\ref{#1}}
\newcommand{\refS}[1]{Section~\ref{#1}}
\newcommand{\refand}[2]{\ref{#1} and~\ref{#2}}

\newcommand\marginal[1]{\marginpar{\raggedright\parindent=0pt\tiny #1}}

\begingroup
  \divide\count255 by 60
  \multiply\count255 by -60
  \advance\count255 by \time
  \ifnum \count255 < 10 \xdef\klockan{\the\count1.0\the\count255}
  \else\xdef\klockan{\the\count1.\the\count255}\fi
\endgroup

\newcommand\nopf{\qed}   

\DeclareMathOperator*{\sumx}{\sum\nolimits^{*}}
\DeclareMathOperator*{\sumxx}{\sum\nolimits^{**}}

\newcommand{\sumxtext}{\sum^{*}}
\newcommand{\sumxxtext}{\sum^{**}}
\newcommand{\sumj}{\sum_{j=0}^{m-1}}
\newcommand{\sumk}{\sum_{k=0}^{m-1}}
\newcommand{\sumki}{\sum_{k=1}^{m-1}}

\newcommand\set[1]{\ensuremath{\{#1\}}}

\newcommand\xpar[1]{(#1)}
\newcommand\bigpar[1]{\bigl(#1\bigr)}
\newcommand\Bigpar[1]{\Bigl(#1\Bigr)}

\def\rompar(#1){\textup(#1\textup)}    
\newcommand\xfrac[2]{#1/#2}

\newcommand\parfrac[2]{\Bigpar{\frac{#1}{#2}}}

\def\xexp(#1){e^{#1}}

\newcommand\ntoo{\ensuremath{{n\to\infty}}}

\newcommand\tti{\ensuremath{{t\to1}}}

\newcommand\iid{i.i.d.\spacefactor=1000}    
\newcommand\ie{i.e.\spacefactor=1000}
\newcommand\eg{e.g.\spacefactor=1000}

\newcommand\cf{cf.\spacefactor=1000}
\newcommand{\as}{a.s.\spacefactor=1000}

\newcommand\ii{\mathrm{i}}

\newcommand{\tend}{\longrightarrow}
\newcommand\dto{\overset{\mathrm{d}}{\tend}}

\newcommand\asto{\overset{\mathrm{a.s.}}{\tend}}
\newcommand\eqd{\overset{\mathrm{d}}{=}}

\newcommand\bbC{\mathbb C}

\newcommand\bbZ{\mathbb Z}

\renewcommand\Re{\operatorname{Re}}
\renewcommand\Im{\operatorname{Im}}

\newcommand\E{\operatorname{\mathbb E{}}}
\renewcommand\P{\operatorname{\mathbb P{}}}
\newcommand\Var{\operatorname{Var}}

\newcommand\Exp{\operatorname{Exp}}
\newcommand\Po{\operatorname{Po}}
\newcommand\Bi{\operatorname{Bi}}

\newcommand\fall[1]{^{\underline{#1}}}

\newcommand\ga{\alpha}
\newcommand\gb{\beta}
\newcommand\gd{\delta}
\newcommand\gD{\Delta}
\newcommand\gf{\varphi}
\newcommand\gam{\gamma}

\newcommand\gl{\lambda}

\newcommand\go{\omega}

\newcommand\gs{\sigma}
\newcommand\gss{\sigma^2}
\newcommand\gS{\Sigma}
\newcommand\eps{\varepsilon}

\newcommand\bZ{\bar Z}

\newcommand\cT{{\mathcal T}}

\newcommand\xhat{\widetilde}

\newcommand\hZ{\xhat Z}
\newcommand\hWm{\xhat W_m}
\newcommand\hZm{\xhat Z_m}

\newcommand\hf{\widehat f}
\newcommand\hg{\widehat g}
\newcommand{\hxm}{\widehat{X}^{(m)}}
\newcommand{\hbxm}{\widehat{\mathbf X}^{(m)}}
\newcommand{\hym}{\widehat{Y}^{(m)}}

\def\[#1]{[\![#1]\!]}

\newcommand\smatrixx[1]{\left(\begin{smallmatrix}#1\end{smallmatrix}\right)}

\newcommand\qq{^{1/2}}
\newcommand\qqm{^{-1/2}}
\newcommand\qm{^{-1}}
\newcommand\qmm{^{-2}}

\renewcommand{\=}{:=}

\newcommand\rv{random variable}
\newcommand\gv{Gaussian variable}

\newcommand\rhs{right hand side}

\newcommand{\GW}{Galton--Watson}
\newcommand{\GWt}{\GW{} tree}
\newcommand{\GWp}{\GW{} process}
\newcommand{\cGWt}{conditioned \GW{} tree}
\newcommand{\Polya}{P\'olya}
\newcommand\gpu{generalized \Polya{} urn}

\newcommand{\xx}[1]{X^{(#1)}}
\newcommand{\xm}{\xx{m}}
\newcommand{\bxx}[1]{\mathbf X^{(#1)}}
\newcommand{\bxm}{\bxx{m}}
\newcommand{\yy}[1]{Y^{(#1)}}
\newcommand{\ym}{\yy{m}}
\newcommand{\byy}[1]{\mathbf Y^{(#1)}}
\newcommand{\bym}{\byy{m}}

\newcommand\etta{\mathbf1} 
\newcommand\gom{\go_m}
\newcommand\vvv{v^{\phantom*}}
\newcommand\phix[1]{\varphi^{(#1)}}
\newcommand\phil{\phix{l}}
\newcommand\gdd{$\gD$-domain}
\newcommand\gda{$\gD$-analytic}

\newcommand{\FT}{Fourier transform}
\newcommand{\DFT}{discrete Fourier transform}

\newcommand{\Cramer}{Cram\'er}

\newcommand{\maple}{\texttt{Maple}}

\newcommand\REM[1]{\texttt{[#1]}\marginal{XXX}}

\newcommand\urladdrx[1]{\urladdr{\def~{\~{}}#1}}

\begin{document}
\title
{Congruence properties of depths in some random trees}

\date{September 21, 2005} 

\author{Svante Janson}
\address{Department of Mathematics, Uppsala University, PO Box 480,
SE-751~06 Uppsala, Sweden}
\email{svante.janson@math.uu.se}
\urladdrx{http://www.math.uu.se/~svante/} 


\begin{abstract} 
Consider a random recusive tree with $n$ vertices.
We show that the number of vertices with even depth is asymptotically
normal as $n\to\infty$. The same is true for the number of vertices
of depth divisible by $m$ for $m=3$, 4 or 5; in all four cases the
variance grows linearly. On the other hand, for $m\ge7$, the number is
not asymptotically normal, and the variance grows faster than linear
in $n$. The case $m=6$ is intermediate: the number is asymptotically
normal but the variance is of order $n\log n$.

This is a simple and striking example of a type of phase transition
that has been observed by other authors in several cases.
We prove, and perhaps explain, this non-intuitive behavious using a
translation to a \gpu.

Similar results hold for a random binary search tree; now the number
of vertices of depth divisible by $m$ is asymptotically normal for
$m\le8$ but not for $m\ge9$, and the variance grows linearly in the
first case both faster in the second. (There is no intermediate case.)

In contrast, we show that for \cGWt{s}, including random labelled
trees and random binary trees, there is no such phase transition:
the number is asymptotically normal for every $m$.
\end{abstract}

\maketitle

\section{Introduction}\label{S:intro}

Given a rooted tree $T$, let $X_j(T)$ be the number of vertices of
depth $j$ (\ie, of distance $j$ from the root).
The sequence $(X_j)_{j=0}^\infty$ is called the \emph{profile} of the tree
and has been studied for various types of random trees by many
authors, see \eg{} 
\cite{AldousII,
ChauvinDJ,
ChauvinKMR,
DrmotaG,
DrmotaH:bi,
DrmotaH:prof,
FuchsHwN,
LouSzT}.

We will here study the congruence class of the depth modulo some given
integer $m\ge2$. Thus, let $\xm_j(T)\=\sum_{k\equiv j \pmod m} X_k(T)$ 
be the number of vertices of depth congruent to $j$ modulo $m$.
For example, $\xx2_0(T)$ is the number of vertices of even depth.
We further let $\bxm(T)$ denote the vector 
$\bigpar{\xm_0(T),\dots,\xm_{m-1}(T)}$.

The purpose of this paper is to study the asymptotic distribution as
\ntoo{} of the random vector $\bxm(T_n)$ for some random trees $T_n$ with
$n$ vertices.
For a random recursive tree (RRT), see \refS{Sresults} for definitions,
we will show (\refT{Trrt})
that $\bxm$ is
asymptotically normal for $m\le6$ but not for $m\ge7$; furthermore,
the variance grows linearly for $m\le5$ but faster for $m\ge6$.
For a random binary search tree (BST), the result is similar
(\refT{Tbst}):
$\bxm$ is
asymptotically normal 
for $m\le8$, but not for $m\ge9$; moreover,
the variance grows linearly for $m\le8$ but faster for larger $m$.
In contrast, for a \cGWt{} (CGWT),
$\bxm$ is asymptotically normal for every $m$  (\refT{Tgwt}).

Note that the typical depths are of the order $\log n$ in a RRT or
BST, but $\sqrt n$ in a CGWT. There is thus more room for smoothing
between the congruence classes in the latter case, which may explain
the asymptotic normality in that case,
but we see no
intuitive explanation for the difference between small and large $m$
for the RRT and BST.

We do not claim that the variables $\xm_j$ have any importance in
applications, but they provide a simple and surprising example of a
type of phase transition that has been observed in several 
similar combinatorial situations.
One well-known such example
is random $m$-ary search trees, see
Chern and Hwang \cite{CHwang}; see also \eg{}
Mahmoud and Pittel \cite{MP}, 
Lew and Mahmoud \cite{LewM}, 
Fill and Kapur \cite{FillK},
Chauvin and Pouyanne \cite{ChauPou}.
Other examples are random quadtrees, see Chern, Fuchs and Hwang
\cite{CFHwang}, and random fragmentation trees, see
Dean and Majumdar \cite{DeanMaj} and 
Janson and Neininger \cite{SJcont}.

We will use a translation into \gpu{s}; for such urns, there is a general
theorem describes the phase transition in terms of the size of the
real part of the second largest eigenvalue, see 
Athreya and Karlin \cite{AK} or Athreya and Ney \cite[\S V.9]{AN};
see also 
Kesten and Stigum \cite{KS:add}
and, for more details,
\cite{SJ154}. The same type of phase transition in a related
(and overlapping) setting is described in \cite{SJ32}.
We hope that our simple example can help to illustrate 
this surprising and non-intuitive phenomenon.
It provides also an illustration of the results in \cite{SJ154} in a
simple concrete situation.

We state the results in \refS{Sresults}. 
Proofs are given in the following three sections, with one type of
random trees in each. The proofs for RRT and BST in Sections
\refand{Srrt}{Sbst} are very similar and uses 
\gpu{s}, while the proof for CGWT in \refS{Sgwt} uses generating
functions
and singularity analysis.
Finally, in \refS{Sosc} we show that the oscillations shown in the
results for RRT and BST except for small $m$ are genuine and not only
an artefact of the proof or of the normalization.

\begin{ack}
  This research was mainly done during the workshop
\emph{Probability theory on trees and analysis of algorithms}
at the Mathematisches Forschungsinstitut Oberwolfach
in August 2004. I thank the participants, in particular Hsien-Kuei
Hwang, for inspiring comments.
\end{ack}

\section{Definitions and results}\label{Sresults}

Let $\etta$ denote the $m$-dimensional vector $(1,\dots,1)$.
(We do not always distinguish between row and column vectors in our notation.)

We will occasionally deal with complex random variables.
A \emph{complex Gaussian variable} is a complex \rv{} $\zeta$ such
that $\Re\zeta$ and $\Im\zeta$ are jointly Gaussian (\ie, normal).
A complex Gaussian variable $\zeta$ is \emph{symmetric} 
if and only if
$\E\zeta=\E\zeta^2=0$; 
its distribution is then determined by the scale factor $\E|\zeta|^2$.
See further \cite[Section I.4]{SJIII}.

The \emph{\DFT} on the group $\bbZ_m=\set{0,\dots,m-1}$
is defined by
\begin{equation}\label{ft}
  \hf(k)=\sumj \go^{kj}f(j),
\end{equation}
where $\go=\gom\=e^{2\pi\ii/m}$.

All limits below are as \ntoo.

\subsection{Random recursive trees}
A \emph{random recursive tree (RRT)} with $n$ vertices 
is a random rooted tree obtained by
starting with 
a single root and then
adding $n-1$ vertices one by one, each time joining the new vertex to a
randomly chosen old vertex;
the random choices are uniform and independent of each other.
If the vertices are labelled $1,2,\dots$, we thus obtain a tree where
the labels increase along each branch as we travel from the root;
the random recursive tree can also be defined as a (uniform) randomly
chosen such labelled tree. 
(The distribution of a random recursive tree differs from
the distribution of a uniform random labelled tree.)
See also the survey \cite{SmytheM}.

We state the main results for RRT in the following theorem, and add
various details in the remarks after it.
\begin{theorem}
  \label{Trrt}
Let $T_n$ be a random recursive tree with $n$ vertices. Then, the
following holds for $\bxm=\bxm(T_n)$, as \ntoo.
  \begin{thmenumerate}
\item
If\/ $2\le m\le 5$, then
\begin{equation}
  \label{trrt1}
n\qqm\Bigpar{\bxm-\frac{n}{m}\mathbf1}\dto N(0,\gS_m)
\end{equation}
for a covariance matrix $\gS_m$ given explicitly in 
\eqref{gs2}--\eqref{gs5}.
\item 
If\/ $ m=6$, then
$$
(n\ln n)\qqm\Bigpar{\bxm-\frac{n}{m}\mathbf1}\dto N(0,\gS_6),
$$
where $\gS_6$ is given explicitly in \eqref{gs6}.
\item
If\/ $ m\ge7$, let $\ga=\cos (2\pi/m)>1/2$ and $\gb=\sin (2\pi/m)$. Then
\begin{equation}\label{trrt3a}
n^{-\ga}
\Bigpar{\bxm- \frac nm \etta}
-\Re\bigpar{ n^{\ii\gb} \hWm} \asto0,
\end{equation}
for some complex random vector 
$\hWm=\bigpar{\hZm e^{-2\pi j\ii/m}}_{j=0}^{m-1}$,
where $\hZm$ is a complex random variable.
In particular, 
along any subsequence such that $\gb\ln n \mod{2\pi}\to\gam$
for some $\gamma\in[0,2\pi]$,
\begin{equation}\label{trrt3b}
n^{-\ga}
\Bigpar{\xm_j- \frac nm }
\dto \Re\bigpar{ e^{\ii(\gam-2\pi j/m)}\hZm},
\end{equation}
jointly in $j=0,\dots,m-1$.
  \end{thmenumerate}
\end{theorem}

\begin{remark}\label{Rrrt}
The covariance matrices $\gS_m$ in (i) and (ii) are given by 
\allowdisplaybreaks
\begin{align}
\Sigma_2 \label{gs2}
&=
\frac1{12}
\left(\!\!\!
\begin{array}{rr}
1 & -1 \\
-1 & 1 
\end{array}
\!\!
\right)
\\
\Sigma_3 
&=
\frac1{18}
\left(\!\!\!
\begin{array}{rrr}
2 & -1 & -1 \\
-1 & 2 & -1 \\
-1 & -1 & 2 
\end{array}
\!\!
\right)
\\
\Sigma_4 
&=
\frac1{48}
\left(\!\!\!
\begin{array}{rrrr}
7 & -1 & -5 & -1 \\
-1 & 7 & -1 & -5 \\
-5 & -1 & 7 & -1 \\
-1 & -5 & -1 & 7
\end{array}
\!\!
\right)
\\
\Sigma_5 \label{gs5}
&=
\frac1{25}
\left(\!\!\!
\begin{array}{rrrrr}
6 & 1 & -4 & -4 & 1 \\
1 & 6 & 1 & -4 & -4 \\
-4 & 1 & 6 & 1 & -4 \\
-4 & -4 & 1 & 6 & 1 \\
1 & -4 & -4 & 1 & 6 
\end{array}
\!\!
\right)
\\
\Sigma_6 \label{gs6}
&=
\frac1{36}
\left(\!\!\!
\begin{array}{rrrrrr}
2 & 1 & -1 & -2 & -1 & 1\\
1 & 2 & 1 & -1 & -2 & -1 \\
-1 & 1 & 2 & 1 & -1 & -2 \\
-2 & -1 & 1 & 2 & 1 & -1 \\
-1 & -2 & -1 & 1 & 2 & 1 \\
1 & -1 & -2 & -1 & 1 & 2 
\end{array}
\!\!
\right)
\end{align}
Note that the matrices $\gS_m$ are circulant, \ie{} invariant under a
cyclic shift of both rows and columns. In particular, $\xm_j$ has the same
asymptotic distribution for every $j$ (when $m\le 6$).

\end{remark}

\begin{remark}
  \label{Rrrt3}
The distribution of $\hZm=\hZ$ in (iii) is determined by the equalities in
distribution
\begin{align}
  Z&\eqd \tfrac m2 W^\go \hZ, \label{hz}
\\
Z&\eqd U^\go(Z+\go Z'), \label{zz}
\end{align}
where $\go=\ga+\gb\ii=e^{2\pi \ii/m}$,
$Z'\eqd Z$, $W\sim\Exp(1)$, $U\sim U(0,1)$, 
and $W$, $\hZ$, $U$, $Z$, $Z'$ are independent,
together with $\E Z=1$.

The moments of $Z$ and $\hZ$ are all finite, and can be obtained
recursively from \eqref{zz} and \eqref{hz} together with $\E Z=1$.
In particular, 
\begin{align*}
\E Z\,\,&=1,   &
\E\hZ&=\frac2{m\Gamma(1+\go)},
\\
\E Z^2\,&=\frac2{2-\go},    &     
  \E \hZ^2&=\frac8{m^2(2-\go)\Gamma(1+2\go)}, 
\\
\E |Z|^2&=\frac{2\ga}{2\ga-1},    &     
  \E |\hZ|^2&=\frac{8\ga}{m^2(2\ga-1)\Gamma(1+2\ga)}, 
\\
\E Z^3&=\frac{6(1+\go)}{(3-\go^2)(2-\go)},    &   
\E \hZ^3&=  
\frac{48(1+\go)}{m^3(3-\go^2)(2-\go)\Gamma(1+3\go)}.
\end{align*}
\end{remark}

\begin{remark}\label{Rrrtmom}
  It follows from \cite[Section 2]{SJ32} that all moments of the
  variables in \refT{Trrt} stay bounded as \ntoo, and thus moment
  convergence holds in (i) and (ii) and for convergent subsequences in
  (iii).
In particular, the variance of each $\xm_j$ is of order $n$ for
  $m\le5$, but larger for $m\ge6$.
\end{remark}

\begin{remark}
The covariance matrix $\gS_m$ has rank $m-1$ when $m\le5$, while
$\gS_6$ has rank 2 only. 
(Full rank is impossible because all row sums are 0, reflecting the
fact that the total number of vertices is non-random.)
\end{remark}

\begin{remark}
  \label{Rrft}
The results for the RRT $T_n$ become simpler if we state them in terms
of the \DFT{}  
$\hbxm=\bigpar{\hxm(k)}_0^{m-1}$ defined by \eqref{ft}. 
(See also \cite{SJ32}, where the proof is based on this \FT.)
Clearly, 
$\hxm(0)=\sum_j \xm_j=|T_n|=n$ is deterministic.

For $2\le m\le 5$, \refT{Trrt} implies (and is equivalent to)
the joint convergence
\begin{equation}\label{rhat}
  n\qqm \hxm(k)\dto V_k,
\qquad k=1,\dots,m-1,
\end{equation}
where $V_k$ are complex (jointly) Gaussian variables such that, for
$k,l\in\set{1,\dots,m-1}$,
using \eqref{sm} below,
\begin{romenumerate}
\item
$V_{m-k}=\overline V_k$, 
\item 
$\E V_k=0$,  
\item
$\E (V_kV_l)=0$ when $l\neq m-k$,
\item
$\displaystyle  
\E|V_k|^2 = \frac1{1-2\Re \gom^k} = \frac1{1-2\cos(2\pi k/m)}$.
\end{romenumerate}
It follows that further:
\begin{romenumerate}
  \setcounter{enumi}{4}
\item
If $k=m/2$ (and thus $\gom^k=-1$),
then $V_k=\overline V_k$ is a symmetric real \gv.
\item
If $k\neq m/2$,
then $\E V_k^2=0$ and thus $V_k$ is a symmetric complex \gv.
\item
The \gv{s} $V_k$, $1\le k\le m/2$, are independent.
\end{romenumerate}
Note that the joint distribution of $V_1,\dots,V_{m-1}$ is determined
by (iv)--(vii) together with (i).

Similarly, for $m=6$, 
\begin{equation*}
  (n \ln n)\qqm \hxm(k)\dto V_k,
\qquad k=1,\dots,5,
\end{equation*}
where now, however, only $V_1$ and $V_{5}=\overline{V_1}$ are
non-zero;
$V_1$ is a symmetric complex \gv{} with $\E|V_1|^2=1$.

Finally, for $m\ge7$, we have
\begin{equation}\label{thomas}
 n^{-\gom} \hxm(1)\asto \frac m2 \hZ_m,
\end{equation}
while $n^{-\gom} \hxm(k)\asto0$ for $k=2,\dots,m-2$.
\end{remark}

\subsection{Binary search trees}

A \emph{Binary search tree (BST)} is constructed from a sequence
$x_1,\dots,x_n$ of distinct real numbers as follows,
see \eg{} \cite[Section 6.2.2]{KnuthIII}.
If $n=0$, the tree is empty. Otherwise, start with a root. (In
computer applications, $x_1$ is stored in the root.)
Then construct recursively two subtrees of the root by the same
procedure applied to two subsequences of $x_1,\dots,x_n$:
the left subtree from the $x_i$ with $x_i<x_1$ and
the right subtree from the $x_i$ with $x_i>x_1$.
The number of vertices in the tree is thus $n$, and each vertex
corresponds to an $x_i$.

A random BST is obtained by this construction applied to a sequence
$x_1,\dots,x_n$ in random order.
(Since only the order properties of $x_1,\dots,x_n$ matter, we can let
them be, for example, either a random permutation of $1,\dots,n$ or
$n$ \iid{} \rv{s} with a common continuous distribution.)

It is easily seen that a random BST can be grown by adding vertices
one by one according to a Markov process similar to the definition of
the RRT:
Given a binary tree with $n$ vertices, there are $n+1$ possible
positions for a new vertex, and we choose one of them at random
(uniformly).

\begin{theorem}
  \label{Tbst}
Let $T_n$ be a random binary search tree with $n$ vertices. Then, the
following holds for $\bxm=\bxm(T_n)$, as \ntoo.
  \begin{thmenumerate}
\item
If\/ $2\le m\le 8$, then
$$
n\qqm\Bigpar{\bxm-\frac{n}{m}\mathbf1}\dto N(0,\gS_m)
$$
for a covariance matrix $\gS_m$ given explicitly in 
\refR{Rbst}.
\item
If\/ $ m\ge9$, let $\ga=2\cos (2\pi/m)-1>1/2$ and $\gb=2\sin (2\pi/m)$. Then
\begin{equation}
n^{-\ga}
\Bigpar{\bxm- \frac nm \etta}
-\Re\bigpar{ n^{\ii\gb} \hWm} \asto0,
\end{equation}
for some complex random vector 
$\hWm=\bigpar{\hZm e^{-2\pi j\ii/m}}_{j=0}^{m-1}$,
where $\hZm$ is a complex random variable.
In particular, 
along any subsequence such that $\gb\ln n \mod{2\pi}\to\gam$
for some $\gamma\in[0,2\pi]$,
\begin{equation}
n^{-\ga}
\Bigpar{\xm_j- \frac nm }
\dto \Re\bigpar{ e^{\ii(\gam-2\pi j/m)}\hZm},
\end{equation}
jointly in $j=0,\dots,m-1$.
  \end{thmenumerate}
\end{theorem}

\begin{remark}
  \label{Rbst}
The covariance matrices $\gS_m$ in (i) 
are circulant and explicitly given by the following first rows:
\begin{align*}
  m=2:&\quad\tfrac1{28}(1,-1),
\\
  m=3:&\quad\tfrac1{45}(2, -1, -1),
\\
m=4:&\quad\tfrac1{336}(17, -3, -11, -3),
\\
m=5:&\quad\tfrac1{275}(16, 1, -9, -9, 1),
\\
m=6:&\quad\tfrac1{1260}(89, 23, -37, -61, -37, 23),
\\
m=7:&\quad\tfrac1{637}(62, 27, -15, -43, -43, -15, 27),
\\
m=8:&\quad\tfrac1{1344}(269, 165, -11, -171, -235, -171, -11, 165).
\end{align*}
\end{remark}

\begin{remark}
  \label{Rbst2}
The distribution of $\hZm=\hZ$ in (ii) is determined by the equalities in
distribution  
\begin{align}
  Z&\eqd \tfrac m2 (2\go-1)W^{2\go-1} \hZ, \label{bhz}
\\
  Z&\eqd \go U^{2\go-1}(Z+Z'), \label{zzbex}
\end{align}
where $\go=\gom=e^{2\pi \ii/m}$,
$Z'\eqd Z$, $W\sim\Exp(1)$, $U\sim U(0,1)$, 
and $W$, $\hZ$, $U$, $Z$, $Z'$ are independent,
together with $\E Z=1$.

Again, the moments of $Z$ and $\hZ$ are all finite, and can be obtained
recursively from \eqref{zzbex} and \eqref{bhz} together with $\E Z=1$.
In particular,
\begin{align*}
\E Z\,\,&=1,   &
\E\hZ&=\frac2{m(2\go-1)\Gamma(2\go)},
\\
\E Z^2\,&=\frac{2\go^2}{4\go-1-2\go^2},    &     
  \E \hZ^2&=\frac4{m^2(2\go-1)^2\Gamma(4\go-1)}\E Z^2, 
\\
\E |Z|^2&=\frac{2}{4\ga-3},    &     
  \E |\hZ|^2&=\frac{4}{m^2(5-4\ga)\Gamma(4\ga-1)}\E|Z|^2, 
\\
\E Z^3&=\frac{6\go^5}{(3\go-1-\go^3)(4\go-1-\go^2)},    
\hskip-7.5pt
&   
\E \hZ^3&=  
\frac{8}{m^3(2\go-1)^3\Gamma(6\go-2)}\E Z^3.
\end{align*}
\end{remark}

\begin{remark}\label{Rbstmom}
The proofs in \cite[Section 2]{SJ32} apply to this situation too
and show
that all moments of the
  variables in \refT{Tbst} stay bounded as \ntoo, and thus moment
  convergence holds in (i)  and for convergent subsequences in
  (ii).
In particular, the variance of each $\xm_j$ is of order $n$ for
  $m\le8$, but larger for $m\ge9$.
\end{remark}

\begin{remark}
  \label{Rbft}
Just as for the RRT case in \refR{Rrft},
the results become simpler if we state them in terms of the \FT{}  
$\hbxm$. Clearly, again
$\hxm(0)=|T_n|=n$.

Further, for $2\le m\le 8$,
\eqref{rhat} holds,
where $V_k$ are complex Gaussian variables
satisfying 
(i)--(iii) and (v)--(vii) in \refR{Rrft} together with
\begin{romenumerate}
\item[(iv$'$)]
$\displaystyle  
\E|V_k|^2 = 
\frac1{3-4\cos(2\pi k/m)}$.
\end{romenumerate}

Finally, for $m\ge9$, 
\begin{equation}\label{ton}
 n^{-(2\gom-1)} \hxm(1)\asto \frac m2 \hZm,
\end{equation}
while $n^{-(2\gom-1)} \hxm(k)\asto0$ for $k=2,\dots,m-2$.
\end{remark}

\subsection{Conditioned \GWt{s}}

A \emph{\cGWt{} (CGWT)} with $n$ vertices is
a random tree obtained as the family tree of a \GWp{}
conditioned on a given total population of $n$. 
(See \eg{} \cite{AldousII,Devroye} for details.)
The \GWp{} is defined using an offspring distribution;
let $\xi$ denotes a random variable with this distribution.
We assume, as usual,
$\E \xi=1$ (the \GWp{} is critical)
and $ 0<\gss=\Var\xi <\infty$. We assume further, for technical reasons, that 
$\E e^{a\xi}<\infty$ for some $a>0$. 
(Equivalently, the probability generating function $\gf(z)\=\E z^\xi$
is analytic in a
disc with radius greater than 1.)

It is well-known \cite{AldousII} that the \cGWt{s}
are the same as the simply generated trees \cite{MM}. 
Many combinatorially interesting random trees are of this type, with
different choices of $\xi$,
for example
labelled trees ($\xi\sim\Po(1)$, $\gss=1$);
ordered (=plane) trees ($\P(\xi=k)=2^{-k-1}$, $\gss=2$);
binary trees ($\xi\sim\Bi(2,1/2)$, $\gss=1/2$);
complete binary trees ($\P(\xi=0)=\P(\xi=2)=1/2$, $\gss=1$).

It has been shown by Drmota and Gittenberger \cite{DrmotaG}
that the profile of a CGWT{} converges, after
normalization, to the local time of a Brownian excursion.
For the congruence classes we have a simpler result. As in many other results,
see \eg{} \cite{AldousII,DrmotaG}, different choices of $\xi$
only affects a scaling factor in the limit result.

\begin{theorem}
  \label{Tgwt}
Let $T_n$ be a random \cGWt{}  with $n$ vertices,
with the offspring distribution given by a random variable $\xi$ such
that
$\E\xi=1$, $\gss\=\Var \xi>0$ and $\E e^{a\xi}<\infty$ for some $a>0$. 
Then, for any fixed $m\ge2$, the
following holds for $\bxm=\bxm(T_n)$ as \ntoo, 
$$
n\qqm\Bigpar{\bxm-\frac{n}{m}\mathbf1}\dto 
N(0,\gS_m)
$$
for a covariance matrix $\gS_m=(\gs_{ij})$ given explicitly by
\begin{equation*}
  \gs_{i,i+k} = \frac{\gss}{12 m^2}\bigpar{m^2-1-6k(m-k)}
\end{equation*}
for $0\le i\le m-1$, $0\le k\le m$ and $i+k$ taken modulo $m$.
\end{theorem}

\begin{remark}
  \label{Rgwt}
The covariance matrices  $\gS_m$ are circulant.
For small $m$, 
they are explicitly given by the following first rows:
\begin{align*}
  m=2:&\quad\tfrac{\gss}{16}(1,-1),
\\
  m=3:&\quad\tfrac{\gss}{27}(2, -1, -1),
\\
  m=4:&\quad\tfrac{\gss}{64}( 5, -1, -3, -1),
\\
m=5:&\quad\tfrac{\gss}{25}(2, 0, -1, -1, 0),
\\
m=6:&\quad\tfrac{\gss}{432}(35, 5, -13, -19, -13, 5).
\end{align*}
\end{remark}

\begin{remark}\label{Rgwtmom}
  It follows from the proof in \refS{Sgwt} that moment convergence
  holds in \refT{Tgwt}.
In particular, the variance of each $\xm_j$ is of order $n$ for
every  $m$.
\end{remark}

\begin{remark}
  Note the curious fact that the asymptotic variance of $\xm_j$,
  $n\gss(1-m\qmm)/12$, is almost independent of $m$.
\end{remark}

\begin{remark}
  \label{Rgt}
For the \FT{}  
$\hbxm$, we again have
$\hxm(0)=|T_n|=n$.
Further, for any $m\ge 2$,
\eqref{rhat} holds,
where $V_k$ are complex Gaussian variables
satisfying 
(i)--(iii) and (v)--(vii) in \refR{Rrft} together with
\begin{romenumerate}
\item[(iv$''$)]
$\displaystyle  
\E|V_k|^2 = 
\frac1{|1-\gom^k|^2}
=\frac1{2-2\cos(2\pi k/m)}$.
\end{romenumerate}
\end{remark}

\begin{remark}
  As a comparison to the results above, suppose that we construct
  $\bxm$ by randomly throwing $n$ balls into $m$ urns $0,\dots,m-1$,
  and counting the number of balls in each urn.
Then $\xm_j\sim\Bi(n,1/m)$, and the central limit theorem shows that 
$n\qqm\bigpar{\bxm-\frac{n}{m}\mathbf1}\dto N(0,\gS_m)$ with
$\gS=(m\qm\gd_{ij}-m\qmm)_{ij}$.
Moreover, the \FT{} $\hbxm$ is a sum of $n$ \iid{} complex \rv{s} and
the central limit theorem shows that 
\eqref{rhat} holds,
where $V_k$ are complex Gaussian variables
satisfying 
(i)--(iii) and (v)--(vii) in \refR{Rrft} together with
\begin{romenumerate}
\item[(iv$'''$)]
$\displaystyle  
\E|V_k|^2 = 1.$
\end{romenumerate}

It follows from the formulas above that the asymptotic variance of
$\xm_j$ (for any $j$) is 
smaller for the depths in the random trees that we consider
than for a random assignment of labels when $m$ is small,
but not when $m$ is large. More precisely, it is smaller for RRT when
$m\le5$ and for BST when $m\le7$; for CGWT, it is smaller when
$m< 12/\gss-1$. (Thus, if $\gss>4$, then the asymptotic variance for
CGWT is always larger than for a random labelling.)
\end{remark}

\section{Random recursive trees}\label{Srrt}

The definition of RRT in \refS{Sresults}
shows immediately that the distribution of
depths modulo $m$ is given by the following \gpu: 
The urn
contains balls with labels $0,\dots,m-1$, representing the depths
modulo $m$ of the vertices.
Start with a single ball with label $0$ in the urn. 
Then, repeatedly,
draw a ball (at random) from the urn, replace it, and if the drawn
ball had label 
$j$, add a new ball with label $j+1 \pmod m$.

This urn was studied briefly in \cite[Example 7.9]{SJ154}, see also
\cite[Example 6.3]{SJ32}. 
We repeat the analysis in \cite{SJ154} with more details.
Using the notation there we have, with the indices in
\set{0,\dots,m-1} and addition taken modulo $m$,
\begin{align*}
\xi_{ij}=\gd_{i+1,j},
&&
a_j=1,
&&
A=(\xi_{ji})_{i,j=0}^{m-1}
=(\gd_{i,j+1})_{i,j=0}^{m-1}.
&
\end{align*}
The matrix $A$ is circulant and corresponds to a convolution operator
in $\ell^2(\bbZ_m)$; hence, $A$ is diagonalized by the characters of
the group $\bbZ_m$. More precisely, let $\go=\gom\=e^{2\pi\ii/m}$.
Then $A$ has left eigenvectors 
$u_j
=(\go^{jk})_{k=0}^{m-1}$
and right eigenvectors $v_j
=(m\qm\go^{-jk})_{k=0}^{m-1}$
with eigenvalues $\go^j$, for $j=0,\dots,m-1$.
Thus, $A$ has $m$ simple eigenvalues, and 
the dominant eigenvalue $\gl_1$ is $\go^0=1$.
(Note that the corresponding eigenvectors $u_0=\etta$ and
$v_0=m\qm\etta$ are denoted 
$u_1$ and $v_1$ in \cite{SJ154}.)
We have chosen the
normalizations such that $(u_j)_j$ and $(v_j)_j$ are
dual bases; moreover $u_0$ and $v_0$ satisfy the normalizations in 
\cite[(2.2)--(2.3)]{SJ154}.

The eigenvalues with
second largest real part are $\gl_2=\go$ and $\gl_3=\overline\go$.
Since $\Re\gl_2=\cos(2\pi/m)<1/2$ when $m\le5$,
$\Re\gl_2=\cos(2\pi/m)=1/2$ when $m=6$,
and
$\Re\gl_2=\cos(2\pi/m)>1/2$ when $m\ge7$,
the trichotomy in \refT{Trrt} follows from
\cite[Theorems 3.22--3.24]{SJ154}.
(The conditions (A1)--(A6) there are easily verified, see
\cite[pp.\ 180--181]{SJ154}.)

More precisely, in Case (i), the convergence to a normal distribution
follows by \cite[Theorems 3.22]{SJ154}.
To find the covariance matrices $\gS_m$, we use 
\cite[Lemma 5.3(iii)]{SJ154}
(or \cite[Lemma 5.3(i) or (ii), Lemma 5.4 and Lemma 5.5]{SJ154})
. 
We have $D=m\qm I$ and thus 
$u_j'Du_k=m\qm u_j\cdot u_k=0$ unless $j\equiv -k\pmod m$,
while
$u_{m-k}'Du_k=u_k^*Du_k=1$.
(We use the notation in \cite{SJ154} that $u\cdot v\=u^tv$, without
complex conjugation.)
Hence,
\begin{equation}\label{sm}
  \gS_m=\sum_{j=1}^{m-1}\frac {1}{1-2\Re \gom^j} \vvv_j v_j^*,
\end{equation}
and a straightforward evaluation yields the matrices \eqref{gs2}--\eqref{gs5}.
(Note that $\gS_m$ is rational also for $m=5$, although $\Re
\go_5=(\sqrt5-1)/4$ is not. This is easily explained by Galois theory.)

The case $m=2$, with $A=\smatrixx{0&1\\1&0}$, is a special case of the
so-called Friedman's urn
\cite{Friedman}
(studied already by Bernstein \cite{Bernstein1}),
and the result follows alternatively directly from
Bernstein \cite{Bernstein1}, \cite{Bernstein2}
or Freedman \cite{Freedman}; see also
\cite[Example 3.27]{SJ154}.

In Case (ii),  $m=6$, we similarly find by
\cite[Lemma 5.3(iv)]{SJ154}
\begin{equation*}
  \gS_6= \vvv_1 v_1^*+\vvv_{5} v_{5}^*=6\qmm\bigpar{2\Re \go_6^{j-k}}_{j,k},
\end{equation*}
which yields \eqref{gs6}.

In Case (iii), \eqref{trrt3a}  follows from
\cite[Theorem 3.24]{SJ154}, 
with $d=0$ and $\gl_2=\gom=\ga+\ii\gb$.
Let us now denote the eigenvectors $u_j$ and $v_j$ belonging to
$\go^j$ by $u_{\go^j}$ and $v_{\go^j}$.
Then \cite[Theorem 3.24]{SJ154} further shows
that $\hWm$ belongs to the linear span $E_\go$ of
$v_\go=(m\qm\go^{-j})_j$; hence 
$\hWm=\bigpar{\hZm \go^{-j}}_{j=0}^{m-1}$ as asserted for some 
complex random $\hZm\=m\qm u_\go\cdot\hWm$. 
By \cite[Theorem 3.26]{SJ154}, \eqref{hz} holds with
$Z=u_\go\cdot W_\go$, and $W_\go$ as in 
\cite[Theorem 3.1]{SJ154}. 

Recall that we start the urn with a single ball with label $0$.
Let, as in 
\cite[Theorem 3.9]{SJ154}, 
$W_{\go,i}$ be the limit random variable 
corresponding to $W_\go$ if we instead start with a single ball of
type $i$. By symmetry, $W_{\go,i}$ is obtained from $W_{\go,0}=W_\go$ 
by a cyclic shift of the components, and thus
$Z_i\=u_\go\cdot W_{\go,i}=\go^i Z$.
By \cite[Theorem 3.9(ii)]{SJ154}, 
\begin{equation}
  \label{sofie}
W_{\go,i}\eqd U^\go(W_{\go,i}+\go W'_{\go,i+1})
\end{equation}
with $W'_{\go,i+1}$ distributed as $W_{\go,i+1}$ and independent of
$W'_{\go,i}$.
Consequently,
\begin{equation*}
  Z_i
\eqd U^\go(Z_i+Z'_{i+1})
\eqd U^\go(Z_i+\go Z'_{i}),
\end{equation*}
with $Z_i$, $Z'_{i+1}$, $Z'_{i}$ independent, and taking $i=0$ we obtain
\eqref{zz}. Moreover, $\E Z=u_{\go 0}=1$ by \cite[Theorem 3.10]{SJ154}.

Conversely, \eqref{zz} implies \eqref{sofie} with $W_{\go,i}=\go^i Z v_\go$,
and thus
\cite[Theorem 3.9(iii)]{SJ154} 
implies that the distribution of $Z$ is
determined by \eqref{zz} and $\E Z$. The distribution of $\hZ$ then is
determined by \eqref{hz}.

For higher moments of $Z$, we take moments in \eqref{zz}. For example,
\begin{equation*}
  \E Z^2 = \E U^{2\go} \E(Z+\go Z')^2 
= (1+2\go)\qm\bigpar{(1+\go^2)\E Z^2+2\go(\E Z)^2},
\end{equation*}
which, using $\E Z=1$, gives $\E Z^2=2/(2-\go)$ after rearrangement.
We leave the corresponding calculations for $\E|Z|^2$ and $\E Z^3$ to
the reader. 
See also \cite[Theorem 3.10]{SJ154}. 

The formulas for moments of $\hZ$ then follows from \eqref{hz}, using
$\E W^z=\Gamma(1+z)$ when $\Re z>-1$;
see also \cite[Theorem 3.26]{SJ154} and \cite[Section 2]{SJ32}.

\section{Binary search trees}\label{Sbst}

To describe the profile of the random BST in terms of an urn model, we
make a simple transformation.
A BST with $n$ vertices has $n+1$ possible positions for a new vertex.
We augment the tree by adding $n+1$ new vertices at these positions;
the $n+1$ new vertices are called \emph{external} and the $n$ original
vertices are called \emph{internal}. Thus every internal vertex has
two children, and every external vertex has none.

Thus, $(X_j)_{j=0}^\infty$ is now the profile of the internal
vertices.
We similarly define $(Y_j)_{j=0}^\infty$ as the profile of the
external vertices, and note that, since every internal vertex has
exactly two children,
\begin{equation*}
  2 X_{j-1}= X_{j}+Y_{j}, \qquad j\ge1.
\end{equation*}
For $j=0$ we instead have $X_0+Y_0=X_0=1$. (Recall that the root has
no parent.)
Passing to congruence classes modulo $m$ we thus have, for every $j$,
\begin{equation} \label{erika}
  2 \xm_{j-1}- \xm_{j}=\ym_{j}-\gd_{j0}.
\end{equation}

The growth of the augmented tree can be described as follows:
Choose an external vertex at random, convert it to an internal vertex
and add two new external vertices as its children.
The distribution of depths for external vertices modulo $m$ is thus
given by a \gpu{} similar to the one in \refS{Srrt}, with the
difference that when we draw a ball with label $j$, we remove it and
add two balls with label $j+1$.
(We start with 2 balls with label 1; alternatively, we start with a
single ball with label 0 and make one more draw.)

The matrix $A$ is now
$(2\gd_{i,j+1}-\gd_{i,j})_{i,j=0}^{m-1}$, again with index addition modulo
$m$.
The eigenvectors are the same $u_j$ and $v_j$ as in \refS{Srrt}, but
the corresponding eigenvalue is now $2\go^j-1$.
In particular, the largest eigenvalue (i.e., the one with largest real
part) is $\gl_1=1$ (as before), and
the second largest are $\gl_2=2\go-1$ and $\overline{\gl_2}$, with
real part
$\Re\gl_2=2\cos(2\pi/m)-1$.

Hence, the condition $\Re \gl_2<1/2$ becomes $\cos(2\pi/m)<3/4$,
which holds for $m\le8$, while for $m\ge9$ we have
$\cos(2\pi/m)>3/4$ and thus $\Re\gl_2>1/2$.

We now obtain, exactly as in \refS{Srrt}, normal convergence
of $\bym=\xpar{\ym_0,\dots,\ym_{m-1}}$
when $m\le8$.
More precisely, by 
\cite[Theorem 3.22 and Lemma 5.3(iii) or Lemmas 5.3(ii) and 5.4]{SJ154},
\eqref{trrt1} holds for $\bym$ with
\begin{equation}\label{smbex}
  \gS_m
=\sum_{j=1}^{m-1}\frac {|2\go^j-1|^2}{3-4\Re \go^j} \vvv_j v_j^*
=\sum_{j=1}^{m-1}\frac {5-4\Re\go^j}{3-4\Re \go^j} \vvv_j v_j^*.
\end{equation}
Similarly,
by \cite[Theorem 3.24]{SJ154}, when $m\ge9$,
\eqref{trrt3a} and \eqref{trrt3b} hold for $\bym$, 
for some $\hZm$ and $\hWm=\bigpar{\hZm \go^{-j}}_{j=0}^{m-1}$. 
Further, by \cite[Theorems 3.26, 3.9 and 3.10]{SJ154}, 
\begin{equation}
    Z\eqd \tfrac m2 W^{2\go-1} \hZm, \label{bexhz}
\end{equation}
where $Z$ satisfies \eqref{zzbex} and $\E Z=1$.

To obtain the results for $\bxm$, we use \eqref{erika}. 
It is convenient to solve this convolution equation by taking the \FT,
which yields, for all $k$,
\begin{equation} \label{herika}
  (2 \go^k-1) \hxm(k)=\hym(k)-1.
\end{equation}
In the case $m\le8$, 
we have in analogy with \refR{Rrft},
$\hym(0)=n+1$ and
the joint convergence
\begin{equation}\label{yrhat}
  n\qqm \hym(k)\dto V_k,
\qquad k=1,\dots,m-1,
\end{equation}
where $V_k$ are complex Gaussian variables satisfying
(i)--(iii) and (v)--(vii) in \refR{Rrft} together with
\begin{romenumerate}
\item[(iv$''''$)]
$\displaystyle  
\E|V_k|^2 
= \frac{|2\gom^k-1|^2}{3-4\Re \gom^k} 
= \frac{5-4\Re\gom^k}{3-4\Re \gom^k} 
= \frac{5-4\cos(2\pi k/m)}{3-4\cos(2\pi k/m)}$.
\end{romenumerate}
It follows immediately from \eqref{herika} and \eqref{yrhat} that
\begin{equation}
  n\qqm \hxm(k)\dto (2\go^k-1)\qm V_k,
\qquad k=1,\dots,m-1,
\end{equation}
which yields the statement in \refR{Rbft} (for $m\le8$), with the
meaning of $V_k$ changed.

Similarly, for $m\ge9$, we have,
with $\gl\=2\go-1$,
in analogy with \eqref{thomas},
$n^{-\gl}\hym(1)\asto \frac m2 \hZ_m$,
and thus
\begin{equation*}
 n^{-\gl} \hxm(1)
= n^{-\gl}\gl\qm \hym(1)+o(1)
\asto \gl\qm\frac m2 \hZ_m,
\end{equation*}
while $n^{-\gl} \hxm(k)\asto0$ for $k=2,\dots,m-2$.
We change the meaning of $\hZ_m$ (replacing $\hZ_m$ by $\gl\hZ_m$) 
and write this as
\eqref{ton}, simultaneously changing \eqref{bexhz} to \eqref{bhz}.

\refT{Tbst} now follows by taking the inverse \FT. When $m\le8$, we
obtain \eqref{trrt1} with
\begin{equation}
  \gS_m=\sum_{j=1}^{m-1}\frac {1}{3-4\Re \gom^j} \vvv_j v_j^*,
\end{equation}
which gives the explicit values in \refR{Rbst}
(with some help of \maple).

\begin{remark}
The covariance matrices in \eqref{smbex} for the case of external vertices 
are circulant and explicitly given by the following first rows:
\begin{align*}
  m=2:&\quad\tfrac1{28}(9, -9),
\\
  m=3:&\quad\tfrac1{45}(14, -7, -7),
\\
m=4:&\quad\tfrac1{336}(97, -27, -43, -27),
\\
m=5:&\quad\tfrac1{275}(76, -9, -29, -29, -9),
\\
m=6:&\quad\tfrac1{1260}(353, 11, -109, -157, -109, 11),
\\
m=7:&\quad\tfrac1{637}(202, 41, -43, -99, -99, -43, 41),
\\
m=8:&\quad\tfrac1{1344}(685, 309, -43, -363, -491, -363, -43, 309).
\end{align*}
\end{remark}

\section{Conditioned \GW{} trees}\label{Sgwt}

For \GW{} trees, we use generating functions and singularity analysis.
See \cite{SJ185} for similar arguments.
Given a tree $T$, we define its profile polynomial by
\begin{equation}\label{sx}
  S(x)=S(x;T)\=\sum_{v\in T} x^{d(v)}=\sum_j x^j X_j(T).
\end{equation}
We will first find the asymptotic distribution of $S(x;T_n)$ for $x$
on the unit circle (excluding the trivial case $x=1$).
Note that $S(x)$ for $|x|=1$ is the \FT{} of the sequence $(X_j)$ as a
function on $\bbZ$.

Letting $\cT$ be a random \GWt,
we define the generating functions, for $k\ge0$,
\begin{equation}\label{fk}
  F_k(t;x_1,\dots,x_k) \= \E\Bigpar{t^{|\cT|}\prod_{i=1}^k S(x_i;\cT)}.
\end{equation}
Here $t$ and $x_1,\dots,x_k$ are complex numbers.
(It is also possible to regard $F_k$ as a formal power series, but we
will need analytic functions.)
We regard $x_1,\dots,x_k$ as fixed and consider 
$F_k$ as a function of $t$.
We consider only $x_i$ with $|x_i|\le1$;
then $|S(x_i;\cT)|\le|\cT|$ and the expectation in \eqref{fk} 
exists at least for $|t|<1$. Thus \eqref{fk}
defines $F_k$ as an analytic function of $t$ in the unit disc $|t|<1$.
We will soon see that it can be continued to a larger domain.

Let $D_0$ be the degree of the root.
If we condition on $D_0=q\ge0$,
then the random tree $\cT$ consists of the root plus $q$ branches
$T_1,\dots,T_q$ that are independent and have the same distribution as
$\cT$.
Further, $|\cT|=1+\sum_1^q|T_j|$ and
$S(x;\cT)=1+\sum_{j=1}^q x S(x;T_j)$, and thus, summing over all
sequences $I_0,\dots,I_q$ of disjoint (possibly empty) subsets of
\set{1,\dots,k} with 
$\bigcup_1^q I_j=\set{1,\dots,k}$,
\begin{equation*}
\prod_{i=1}^k S(x_i;\cT)
=\sum_{I_0,\dots,I_q} \prod_{j=1}^q
\prod_{i\in I_j} x_i S(x_i;T_j).
\end{equation*}
Consequently, 
\begin{equation*}
    \begin{split}
\E\Bigpar{t^{|\cT|}\prod_{i=1}^k S(x_i;\cT)\Bigm| D_0=q}
&=
\E \sum_{I_0,\dots,I_q} t  
\prod_{j=1}^q \Bigpar{ t^{|T_j|}
 \prod_{i\in I_j} x_i S(x_i;T_j)}
\\
&=
t \sum_{I_0,\dots,I_q}   
\Bigpar{ \prod_{i\notin I_0} x_i}
\prod_{j=1}^q F_{|I_j|}(t;\set{x_i,\,i\in I_j}).
  \end{split}
\end{equation*} 
The terms in the latter sum do not depend on the order of
$I_1,\dots,I_q$. Thus, if $\sumxtext_{I_0,\dots,I_l}$ denotes the sum over
such sequences $I_0,\dots,I_l$ with $I_1,\dots,I_l$ non-empty and in,
say, lexicographic order, then, with $q\fall l\=q(q-1)\dotsm(q-l+1)$,
\begin{multline}\label{fk-}
\E\Bigpar{t^{|\cT|}\prod_{i=1}^k S(x_i;\cT)\Bigm| D_0=q}
\\
=  
t \sum_{l=0}^k \sumx_{I_0,\dots,I_l}  q\fall l
\Bigpar{ \prod_{i\notin I_0} x_i}
\prod_{j=1}^l F_{|I_j|}(t;\set{x_i,\,i\in I_j})
F_0(t)^{q-l}.
\end{multline}
Now take the expectation, \ie{} multiply by $\P(D_0=q)$ and sum over
$q$. We have, for $|z|<1$ at least, since
$D_0\eqd\xi$ and thus $\E z^{D_0}=\gf(z)$,
$$
\sum_{q\ge l} q\fall l z^{q-l} \P(D_0=q)
=\E\bigpar{D_0\fall l z^{D_0-l}}
=\phil(z),
$$
and thus \eqref{fk-} yields
\begin{equation*}
F_k(t;x_1,\dots,x_k)
=
t \sum_{l=0}^k \sumx_{I_0,\dots,I_l}  \phil(F_0(t))
\Bigpar{ \prod_{i\notin I_0} x_i}
\prod_{j=1}^l F_{|I_j|}(t;\set{x_i,\,i\in I_j}).
\end{equation*}
In particular, $k=0$ yields the well-known formula
$F_0(t)=t\gf(F_0(t))$.
The next two cases are
\begin{align*}
  F_1(t;x)=t\gf(F_0(t))& + t\gf'(F_0(t)) xF_1(t;x),
\\
  F_2(t;x,y)=t\gf(F_0(t)) 
&+ t\gf'(F_0(t))\bigpar{ xF_1(t;x)+yF_1(t;y)+xyF_2(t;x,y)}
\\&
+ t\gf''(F_0(t)) xyF_1(t;x)F_1(t;y).
\end{align*}
We thus have, recalling $t\gf(F_0(t))=F_0(t)$,
\begin{align}
  F_1(t;x)&= \frac{F_0(t)}{1-xt\gf'(F_0(t))},
\label{f1}
\\
  F_2(t;x,y)&=
\bigpar{1-xyt\gf'(F_0(t))}\qm
\Bigpar{
F_0(t) 
+ t\gf'(F_0(t))\bigpar{ xF_1(t;x)+yF_1(t;y)}
\notag\\&\hskip14em
+ t\gf''(F_0(t)) xyF_1(t;x)F_1(t;y)},
\label{f2}
\end{align}
and in general
\begin{multline}
  \label{fk1}
F_k(t;x_1,\dots,x_k)
=
\Bigpar{1-t\Bigpar{\prod_{i=1}^k x_i}\gf'(F_0(t))}\qm
\\\times
t \sum_{l=0}^k \sumxx_{I_0,\dots,I_l}  \phil(F_0(t))
\Bigpar{ \prod_{i\notin I_0} x_i}
\prod_{j=1}^l F_{|I_j|}(t;\set{x_i,\,i\in I_j}),
\end{multline}
where $\sumxxtext$ means $\sumxtext$ with the single term with $l=1$ and
$|I_1|=k$ omitted.
This gives recursively an explicit formula for each 
$F_k(t;x_1,\dots,x_k)$ as a rational function of $x_1,\dots,x_k$
and $t\phil(F_0(t))$, $0\le l\le k$.

We say, see \cite[Chapter VI]{FS}, 
that a \emph{\gdd}
is a domain of the type \set{z:|z|<1+\eps,\, |\arg(z-1)|>\pi/2-\gd}
for some (small) positive $\eps$ and $\gd$, and that a function is
\emph{\gda} if it is analytic in some \gdd, or can be extended to
such a function.

Suppose for simplicity in the sequel
that $\xi$ is aperiodic. (The periodic case is
similar with standard modifications as in \cite[Chapter VI.7]{FS}; we
omit the details.) Then, by a standard result in singularity analysis,
see \eg{} \cite[Proposition VI.1]{FS},
$F_0(t)$ is \gda, with 
\begin{equation}\label{f0}
 F_0(t)=1-\sqrt{2/\gss}(1-t)\qq+O(1-t)\qquad \text{as \tti}. 
\end{equation}
(Here and below, we consider only $t$ in a suitable \gdd.)
It follows that in a (possibly smaller) \gdd, $|F_0(t)|<1$ and hence 
$|\gf'(F_0(t))|<\gf'(1)=1$. 
Further, 
\begin{equation}
  \label{emma}
  \begin{split}
\gf'(F_0(t))
&=1+\gf''(1)\bigpar{F_0(t)-1}+O\bigpar{F_0(t)-1}^2
\\&
=1-\sqrt2\gs(1-t)\qq+O(1-t),	
  \end{split}
\end{equation}
and it follows easily that
$|t\gf'(F_0(t))|<1$ in a \gdd.
Hence, for every fixed $x$ with $|x|\le1$,
$\bigpar{1-xt\gf'(F_0(t))}\qm$
is \gda, with 
\begin{equation}\label{jesper}
\bigpar{1-xt\gf'(F_0(t))}\qm=
\begin{cases}
2\qqm\gs\qm(1-t)\qqm+O(1),& \text{if $x=1$},
\\
O(1),& \text{if $x\neq1$}.
\end{cases}
\end{equation}

\begin{lemma}
  \label{Lfk}
For $k\ge1$ and any complex $x_1,\dots,x_k$ with $|x_i|\le1$ and $x_i\neq1$,
$F_k(t;x_1,\dots,x_k)$ is \gda{} (as a function of $t$), with 
\begin{equation}
  \label{fka}
F_k(t;x_1,\dots,x_k)
=a_k(x_1,\dots,x_k)(1-t)^{-(k-1)/2} + O\bigpar{(1-t)^{-(k-2)/2}},
\end{equation}
where $a_1(x)=1/(1-x)$, and for $k\ge2$,
$a_k(x_1,\dots,x_k)=0$ if $x_1\dotsm x_k\neq1$, 
while if $x_1\dotsm x_k=1$, 
\begin{equation}\label{ak}
a_k(x_1,\dots,x_k)
=2^{-3/2}\gs\!\sum_{\emptyset\subsetneq I \subsetneq \set{1,\dots,k}}
a_{|I|}(x_i:i\in I)a_{k-|I|}(x_i:i\notin I).
\end{equation}
In particular, when $|x|=1$,
$
a_2(x,\bar x)=2\qqm\gs|1-x|\qmm.
$
\end{lemma}

\begin{proof}
  The $\gD$-analyticity follows by \eqref{fk1} and induction, using
  the results just shown.

For $k=1$, the expansion \eqref{fka} follows from \eqref{f1}, \eqref{f0} 
and \eqref{emma}.

Similarly,  \eqref{fka} for $k=2$ follows from
\eqref{f2} together with \eqref{f0}, \eqref{emma}
and \eqref{jesper}; this also yields 
$a_2(x,y)=0$
if $xy\neq1$ and for $xy=1$, using also $\gf''(1)=\gss$,
\begin{equation*}
  a_2(x,y)
=2\qqm\gs\qm\bigpar{F_0(1)+xa_1(x)+ya_1(y)+\gss xya_1(x)a_1(y)},
\end{equation*}
which equals $2\qqm\gs|1-x|\qmm$
because now $|x|=|y|=1$ and $y=\bar x$ and thus 
$$
xa_1(x)+ya_1(y)=\frac{x}{1-x}+\frac{y}{1-y}=2\Re\Bigpar{\frac{x}{1-x}}=-1.
$$

For $k\ge3$ we argue similarly. By induction, all terms in the sum in
\eqref{fk1} are $O\bigpar{(1-t)^{-(k-2)/2}}$. 
The result when
$x_1\dotsm x_k\neq1$
follows immediately by \eqref{jesper}.

Assume now $x_1\dotsm x_k=1$. If $|I_1|=k-1$, then
$I_1=\set{1,\dots,k}\setminus x_p$ for some $p$, and thus
$\prod_{i\in I_1} x_i=x_p\qm\neq1$; hence by induction
$F_{k-1}(t;\set{x_i:i\in I_1})=O\bigpar{(1-t)^{-(k-3)/2}}$. 
The leading terms in \eqref{fk1} are thus those with $l=2$ and
$I_0=\emptyset$, $I_2=\set{1,\dots,k}\setminus I_1$, which proves the
claim, including \eqref{ak}, 
by another application of \eqref{jesper}.
(Note the factor 1/2 because we assume that $I_1$ and $I_2$ are in order,
but not necessarily $I$ and its complement.) 
\end{proof}

We next solve the recursion \eqref{ak}.
Let $N(x_1,\dots,x_k)$ be the number of pairings of $x_1,\dots,x_k$
into $k/2$ pairs of the type \set{x,\bar x}. 
(Thus $N(x_1,\dots,x_k)=0$ if $k$ is odd.)

\begin{lemma}
  \label{Lak} Let $k\ge2$.
Suppose that $x_1,\dots,x_k\in\set{x\in\bbC:|x|=1 \text{ but }x\neq1}$.
Then
\begin{equation*}
  a_k(x_1,\dots,x_k)=\frac{\Gamma\bigpar{\xfrac{(k-1)}2}}{\sqrt{2\pi}\gs}
  N(x_1,\dots,x_k)\prod_{i=1}^k \frac{\gs}{1-x_i}.
\end{equation*}
\end{lemma}

\begin{proof}
  For $k=2$, the result follows directly from \refL{Lfk}.

For $k\ge3$, we use induction.
The result is trivial if $x_1\dotsm x_k\neq1$.
Hence, we assume $x_1\dotsm x_k=1$ and use \eqref{ak}.
First, note that it suffices to consider $I$ with $|I|$ even.
In fact, if $|I|$ is odd and $|I|\ge3$, then 
$N_{|I|}(x_i:i\in I)=0$, and thus
$a_{|I|}(x_i:i\in I)=0$ by induction.
If $|I|=1$, then $I=\set{x_p}$ for some $p$. Since 
$\prod_{i\notin I_1} x_i=x_p\qm\neq1$, we have 
$N_{k-|I|}(x_i:i\notin I)=0$, and thus by induction
$a_{k-|I|}(x_i:i\notin I)=0$.
Similarly, we may assume that $k-|I|$ is even.

The result thus holds when $k$ is odd (with 
$a_k(x_1,\dots,x_k)=0$).

Now let $k$ be even, $k=2l$ with $l\ge2$.
We use induction on the \rhs{} of \eqref{ak}.
Say that a pairing of $x_1,\dots,x_k$ is \emph{good} if each pair
consist of two conjugate numbers. 
Note that
$N_{|I|}(x_i:i\in I)N_{k-|I|}(x_i:i\notin I)$ 
equals the number of pairs $(\ga,\gb)$ where 
$\ga$ is a good pairing of $(x_i:i\in I)$
and $\gb$ is a good pairing of $(x_i:i\notin I)$.
Each such pair $(\ga,\gb)$ defines a good pairing of $x_1,\dots,x_k$,
and conversely, each good pairing of $x_1,\dots,x_k$ splits into good
pairings $\ga$ and $\gb$ in $2^l-2$ ways: 
$\binom lj$ ways with
$|I|=2j$ for each $j=1,\dots,l-1$.
Consequently, \eqref{ak} yields,
\begin{multline*}
a_k(x_1,\dots,x_k)
=2^{-3/2}\gs\sum_{j=1}^{l-1} \binom lj N(x_1,\dots,x_k)
\\\times
\frac{\Gamma\xpar{j-1/2}\Gamma(l-j-1/2)}{{2\pi}\gss}
\prod_{i=1}^k \frac{\gs}{1-x_i}.
\end{multline*}
To complete the induction step, it is now sufficient to verify
\begin{equation*}
\sum_{j=1}^{l-1} \binom lj \Gamma\xpar{j-1/2}\Gamma(l-j-1/2)
=4\sqrt\pi \Gamma(l-1/2),
\end{equation*}
for $l\ge2$. This is an immediate consequence of the binomial
convolution
\begin{equation*}
  \begin{split}
\sum_{j=0}^{l} \binom lj \Gamma\xpar{j-1/2}\Gamma(l-j-1/2)
=
l!\sum_{j=0}^{l} \frac{\Gamma\xpar{j-1/2}}{j!}
\frac{\Gamma(l-j-1/2)}{(l-j)!}
\\
=
l!\sum_{j=0}^{l} (-1)^l \Gamma(-1/2)^2\binom{1/2}j \binom{1/2}{l-j}
=
l!\,(-1)^l \Gamma(-1/2)^2\binom{1}l=0,	
  \end{split}
\end{equation*}
when $l\ge2$,
since the terms with $j=0$ and $j=l$ both are 
$\Gamma(-1/2)\Gamma(l-1/2)=-2\sqrt\pi\Gamma(l-1/2)$.
\end{proof}

\begin{lemma}
  \label{Lsk}
Let $k\ge1$.
Suppose that $x_1,\dots,x_k\in\set{x\in\bbC:|x|=1 \text{ but }x\neq1}$.
Then
\begin{equation*}
\E \bigpar{S(x_1;T_n)\dotsm S(x_k;T_n)}
=   N(x_1,\dots,x_k)\Bigpar{\prod_{i=1}^k \frac{\gs}{1-x_i}} n^{k/2} 
+ O\bigpar{n^{(k-1)/2}}.
\end{equation*}
\end{lemma}

\begin{proof}
  By \eqref{f1},
  \begin{equation*}
\E \bigpar{S(x_1;T_n)\dotsm S(x_k;T_n)}
= \frac{[t^n]F_k(t;x_1,\dots,x_k)}{\P(|\cT|=n)}	
= \frac{[t^n]F_k(t;x_1,\dots,x_k)}{[t^n]F_0(t)}.
  \end{equation*}
Standard singularity analysis
see \eg{} \cite[Chapter VI]{FS},
and \eqref{fka} yield, for $k\ge2$,
\begin{equation*}
[t^n]F_k(t;x_1,\dots,x_k)
=\frac{a_k(x_1,\dots,x_k)}{ \Gamma\bigpar{(k-1)/2}} n^{(k-3)/2}
+O\bigpar{n^{(k-4)/2}},
\end{equation*}
while for $k=1$ we get
\begin{equation*}
[t^n]F_1(t;x)
=O\bigpar{n^{-3/2}}.
\end{equation*}
Similarly, \eqref{f0} yields, as is well-known, 
\begin{equation*}
[t^n]F_0(t)
=\frac1{\sqrt{2\pi}\gs}
n^{-3/2}
+O\bigpar{n^{-5/2}}.
\end{equation*}
The result follows from these formulas and \refL{Lak}.
\end{proof}

We can now identify the asymptotic moments and thus the asymptotic
distribution of $S(x;T_n)$.

\begin{theorem}\label{Tsx}
Let $U(x)$ be a family of complex Gaussian \rv{s}, defined for $|x|=1$
but $x\neq1$, such that
\begin{romenumerate}
\item 
$U(x)$ is symmetric complex Gaussian when $\Im x\neq0$,
\newline
with $\E|U(x)|^2=1/|1-x|^2$;
\item 
$U(x)$ is symmetric real Gaussian when $\Im x=0$
(\ie, when $x=-1$),
with $\E|U(x)|^2=1/|1-x|^2$;
\item 
$U(\bar x)=\overline{U(x)}$;
\item
the variables $U(x)$, $\Im x\ge0$, are independent.
\end{romenumerate}
Then, for the CGWT, $S(x;T_n)/\sqrt n\dto \gs U(x)$, jointly for all such $x$.
\end{theorem}
\begin{proof}
Note that the assumptions imply that $\E U(x)U(y)=1/|1-x|^2$ if $xy=1$
(and thus $y=\bar x$), but $\E U(x)U(y)=0$ otherwise.
By \refL{Lsk} and the formula
\cite[Theorem 1.28]{SJIII} for joint moments of Gaussian variables
(known as Wick's theorem),
\begin{equation*}
n^{-k/2}  \E \bigpar{S(x_1;T_n)\dotsm S(x_k;T_n)}
\to \gs^k \E \bigpar{U(x_1)\dotsm U(x_k)}.
\end{equation*}
Replacing one or several $x_i$ by their conjugates, we see that 
the same holds if we replace some $S$ and $U$ by their conjugates.
Hence the result holds by the method of moments (applied to the real
and imaginary parts).
\end{proof}

\begin{proof}[Proof of \refT{Tgwt}]
As remarked above, for any tree $T$,
the numbers $S(x;T)$ for $|x|=1$ form the \FT{} of the sequence $(X_j(T))$.
It follows that the discrete \FT{} $\hbxm$ of $\bxm(T)$ equals
the vector $(S(\go^k;T))_{k=0}^{m-1}$, 
where $\go=\gom$,
\cf{} \eqref{ft} and
\eqref{sx}.
Hence, by Fourier inversion,
\begin{equation*}
  X_j(T)=\frac1m \sumk \go^{-jk} S(\go^k;T).
\end{equation*}
\refT{Tsx} thus implies, since trivially $S(1;T_n)=n$,
\begin{equation*}
 n\qqm \Bigpar{X_j(T_n)-\frac nm}
\dto
Z_j\=
\frac{\gs}m
\sumki \go^{-jk} U(\go^k),
\end{equation*}
jointly for all $j$.
Since the variables $U(\go^k)$ are jointly (complex) Gaussian,
the variables $Z_j$ are too; moreover, each $Z_j$ is real.
Clearly, $\E Z_j=0$, and the covariance matrix is given by
\begin{equation}
  \label{magnus}
  \begin{split}
\E(Z_iZ_j)=\E(Z_i\overline{Z_j})
=\frac{\gss}{m^2} \sumki \go^{(j-i)k} \E |U(\go^k)|^2
\\
=\frac{\gss}{m^2} \sumki \go^{(j-i)k} |1-\go^k|\qmm.
  \end{split}
\end{equation}

To evaluate this sum, define a function $f$ on 
$\bbZ_m=\set{0,\dots,m-1}$ by $f(j)=j-(m-1)/2$.
Then its \FT{} is
\begin{equation*}
  \hf(k)=\sumj\bigpar{j-\frac{m-1}2}\go^{jk}
=\frac{m}{\go^k-1},
\qquad
k\neq0,
\end{equation*}
while $\hf(0)=0$.

Let further $g:=f*\check f$ on $\bbZ_m$, \ie{}
$g(j)=\sum_{i\in\bbZ_m} f(j+i)f(i)$.
Then
$\hg(k)=|\hf(k)|^2=m^2|\go^k-1|\qmm$, $k\neq0$, and thus, by Fourier
inversion again,
\begin{equation*}
  \sumki \go^{(j-i)k} |1-\go^k|\qmm
=
m\qmm  \sumk \go^{(j-i)k} \hg(k) = m\qm g(i-j).
\end{equation*}
Hence, by \eqref{magnus}, $\E(Z_iZ_j)=\gss m^{-3} g(i-j)$.
It remains only to evaluate $g$. For $0\le j<m$,
\begin{equation*}
  \begin{split}
  g(j)
&
=\sum_{i=0}^{m-1-j} (j+i)i + \sum_{i=m-j}^{m-1} (j+i-m)i -
 m\parfrac{m-1}2^2
\\&
=m\frac{m^2-1}{12}-m\frac{j(m-j)}{2}.	
  \end{split}
\end{equation*}
\vskip-\baselineskip
\end{proof}

\section{Oscillations}\label{Sosc}

In Case (iii) of \refT{Trrt} ($m\ge7$), we do not have convergence in
distribution:
the sequence of random vectors $n^{-\ga}\bigpar{\bxm- \frac nm \etta}$, 
$n=1,2,\dots$, is tight and thus suitable subsequences converge as is
shown explicitly in \eqref{trrt3b}, but different subsequences may have
different limits, and thus
$n^{-\ga}\bigpar{\xm_0- \frac nm \etta}$, for example, does not have a
limiting distribution.
Indeed, suppose that  
$n^{-\ga}\bigpar{\xm_0- \frac nm \etta}\dto V$, say.
Then, by \eqref{trrt3b}, 
$\Re\xpar{ e^{\ii\gam}\hZm}\eqd V$ for every $\gamma\in[0,2\pi]$.
In particular, 
$\Re\xpar{ e^{\ii\gam}\E\hZm}
=\E\Re\xpar{ e^{\ii\gam}\hZm}
=\E V$ is independent of $\gamma\in[0,2\pi]$,
which is a contradicion because $\E\hZm\neq0$ by \refR{Rrrt3}.

Nevertheless, it is conceivable (although implausible)
that $\xm_0$ has a limiting
distribution if we choose the norming constants carefully,
\ie{} that $a_n(\xm_0-b_n)\dto V$ for some non-degenerate $V$
and real constants $a_n>0$ and $b_n$.
It then would follow from \eqref{trrt3b} that
$\Re\xpar{ e^{\ii\gam}\hZm}\eqd a_\gam V+b_\gam$ 
for every $\gamma\in[0,2\pi]$ and some real constants $a_\gam\ge0$ and
$b_\gam$, see \eg{} \cite[Section 9.2]{Gut}.
In other words, $\Re\xpar{ e^{\ii\gam}\hZm}$ would have 
a distribution of the same type
for every $\gamma$, except when it is degenerate.

To rule this out, we show a general result.

\begin{proposition}\label{P1}
Let $Z$ be a complex random variable such that
$\E|Z|^3<\infty$. Suppose that there exists a random variable $V$ and, 
for every $\gamma\in[0,2\pi]$,
some real constants $a_\gam\ge0$ and $b_\gam$ such that
$\Re\bigpar{ e^{\ii\gam}Z}\eqd a_\gam V+b_\gam$.
Then either
\begin{romenumerate}
  \item
$Z\eqd aW+b$ for some real random variable $W$ and some complex
  constants $a$ and $b$, 
and thus $|\E(Z-\E Z)^2|=\E|Z-\E Z|^2$;
or
\item
$\E (Z-\E Z)^3= \E (Z-\E Z)^2\overline{(Z-\E Z)}=0$.
\end{romenumerate}
\end{proposition}


\begin{proof}
If all $a_\gam=0$, then
$\Re\bigpar{e^{\ii\gam}(Z-\E Z)}=0$ \as{} for every $\gam$, and it
follows by the \Cramer--Wold device that $Z=\E Z$ \as, a special case
of (i).

Thus assume that some $a_\gam>0$. Then $\E |V|^3<\infty$.
By replacing $Z$ by $Z-\E Z$ and $V$ by $V-\E V$ (changing $b_\gam$
accordingly),  we may assume that $\E Z=\E V=0$, and thus $b_\gam=0$.
If $V=0$ \as, then, by the \Cramer--Wold device again, $Z=0$ \as,
and (i) holds. Assume thus $\E V^2>0$; rescaling $V$ we may assume $\E V^2=1$.
Define
\begin{align*}
X_\gam&\=\Re\bigpar{ e^{\ii\gam}Z}\eqd a_\gam V,
\\
Z_\gam&\=2 e^{\ii\gam} X_\gam
=e^{2\ii\gam}Z+\bZ.
\end{align*}
We have
\begin{align*}
\E Z_\gam^2
&=4 e^{2\ii\gam} \E X_\gam^2
=4 e^{2\ii\gam}a_\gam^2  ,
\\
\E Z_\gam^3
&=8 e^{3\ii\gam} \E X_\gam^3
=8e^{3\ii\gam}a_\gam^3 \E V^3,
\end{align*}
and thus
\begin{equation}
  \label{manne}
\bigpar{\E Z_\gam^3}^2
=
(\E V^3)^2 \bigpar{\E Z_\gam^2}^3.
\end{equation}
On the other hand,
\begin{align}
\E Z_\gam^2 
&=   
\E\bigpar{e^{2\ii\gam}Z+\bZ}^2
=
e^{4\ii\gam}\E Z^2
+2e^{2\ii\gam}\E|Z|^2
+\E\bZ^2,   \label{q2}
\\
\E Z_\gam^3
&=   
e^{6\ii\gam}\E Z^3
+3e^{4\ii\gam}\E (Z^2\bZ)
+3e^{2\ii\gam}\E(Z\bZ^2)
+\E\bZ^3. \label{q3}
\end{align}
Hence, $\E Z_\gam^2$ and $\E Z_\gam^3$ are polynomials 
$p_2\bigpar{e^{2\ii\gam}}$ and $p_3\bigpar{e^{2\ii\gam}}$
in
$e^{2\ii\gam}$ of degrees at most 2 and 3. 
By \eqref{manne},
\begin{equation}
  \label{bo}
p_3(z)^2=(\E V^3)^2 p_2(z)^3
\end{equation}
for every $z$ with $|z|=1$, and thus
for every complex $z$.

If $\E V^3=0$, then $\E Z_\gam^3=0$ by \eqref{manne} and thus
(ii) holds by \eqref{q3}.

Suppose now $\E V^3\neq0$. If $\E Z^2=0$, then $p_2(z)$ has degree 1
or 0 by \eqref{q2}. Degree 1 is impossible by \eqref{bo},
and thus $\E|Z|^2=0$ by \eqref{q2} so $Z=0$ \as, and both (i) and (ii)
hold.

Finally, if $\E V^3\neq0$ and $\E Z^2\neq0$, then \eqref{bo} implies
that $p_2$ has a double root, so its discriminant
$(\E|Z|^2)^2-\E Z^2\E\bZ^2=0$, and $|\E Z^2|=\E|Z|^2$. 
It follows that the argument of $Z^2$
is constant \as, and thus (i) holds.
\end{proof}

Returning to \refT{Trrt}(iii), we can use the moments computed in
\refR{Rrrt3}.
The proposition shows that there really are
oscillations, even with different normalizations, as soon as
$|\E(\hZm-\E \hZm)^2|<\E|\hZm-\E \hZm|^2$ and
$\E (\hZm-\E \hZm)^3\neq0$.

It should be possible to verify this for all $m\ge7$, perhaps using
asymptotical expansions for large $m$, but for simplicity we have
resorted to numerical verification (by \maple) for $m\le100$.
We have also done the same for BST and \refT{Tbst}(ii), using the
moments given in \refR{Rbst2}.
We thus conclude the following result, showing that at least for these
$m$, there are genuine oscillations.

\begin{theorem}
  For RRT, at least for $7\le m\le 100$, there are oscillations in
  \refT{Trrt}(iii);
$a_n(\xm_0-b_n)$ does not have a non-degenerate limit distribution for
  any sequence of norming constants $a_n\ge0$ and $b_n$.
The same holds for BST, at least for $9\le m\le 100$.
\nopf
\end{theorem}

\begin{remark}
The fact that $\E (\hZm-\E \hZm)^3\neq0$ also implies by \eqref{q3}
that the subsequence limit 
$X_\gam\=\Re\xpar{ e^{\ii\gam}Z}$ has a non-zero third central
moment, except for at most 6 values of $\gamma\in[0,2\pi)$; in
  particular, $X_\gam$ is not normal except possibly for a few
  exceptional $\gam$.
Presumably, these too could be eliminated by considering fourth or
fifth moments, but we have not pursued that.
  \end{remark}

\newcommand\AAP{\emph{Adv. Appl. Probab.} }
\newcommand\JAP{\emph{J. Appl. Probab.} }
\newcommand\AMS{Amer. Math. Soc.}
\newcommand\JAMS{\emph{J. \AMS} }
\newcommand\MAMS{\emph{Memoirs \AMS} }
\newcommand\PAMS{\emph{Proc. \AMS} }
\newcommand\TAMS{\emph{Trans. \AMS} }
\newcommand\AnnMS{\emph{Ann. Math. Statist.} }
\newcommand\AnnPr{\emph{Ann. Probab.} }
\newcommand\CPC{\emph{Combin. Probab. Comput.} }
\newcommand\JMAA{\emph{J. Math. Anal. Appl.} }
\newcommand\RSA{\emph{Random Struct. Alg.} }
\newcommand\SPA{\jour{Stochastic Process. Appl.} } 
\newcommand\ZW{\emph{Z. Wahrsch. Verw. Gebiete} }
\newcommand\DMTCS{\jour{Discr. Math. Theor. Comput. Sci.} }

\newcommand\Springer{Springer}
\newcommand\Wiley{Wiley}

\newcommand\vol{\textbf}
\newcommand\jour{\emph}
\newcommand\book{\emph}
\newcommand\inbook{\emph}
\def\no#1#2,{\unskip#2, no. #1,} 

\newcommand\webcite[1]{\hfil\penalty0\texttt{\def~{\~{}}#1}\hfill\hfill}
\newcommand\webcitesvante{\webcite{http://www.math.uu.se/\~{}svante/papers/}}
\newcommand\arxiv[1]{arXiv:#1}

\def\nobibitem#1\par{}

\end{document}